\documentclass[12pt,a4paper]{amsart}
\usepackage{amssymb}
\usepackage{mathrsfs}

\headsep=1cm \numberwithin{equation}{section}
\newtheorem{thm}{Theorem}[section]
\newtheorem{cor}[thm]{Corollary}
\newtheorem{lem}[thm]{Lemma}
\newtheorem{prop}[thm]{Proposition}
\theoremstyle{definition}

\theoremstyle{remark}

\numberwithin{equation}{section}
\newcommand\Att{\operatorname{Att}}
\newcommand\Max{\operatorname{Max}}
\newcommand\Supp{\operatorname{Supp}}
\newcommand\Ass{\operatorname{Ass}}

\newcommand\Ann{\operatorname{Ann}}
\newcommand\Spec{\operatorname{Spec}}
\newcommand\Rad{\operatorname{Rad}}
\newcommand\Hom{\operatorname{Hom}}
\newcommand\Ext{\operatorname{Ext}}

\begin{document}
\title [Modules cofinite and weakly cofinite]{Modules cofinite and weakly cofinite  with respect to an ideal}
\author[K. Bahmanpour, R. Naghipour  and M. Sedghi ]{Kamal Bahmanpour, Reza Naghipour$^{*,\dag}$  and Monireh Sedghi}
\address{Department of Mathematics,  Faculty of Mathematical
Sciences, University of Mohaghegh Ardabili, 56199-11367,  Ardabil, Iran; and School of Mathematics, Institute for Research in Fundamental
Sciences (IPM), P.O. Box. 19395-5746, Tehran, Iran.}
\email{bahmanpour.k@gmail.com}
\address{Department of Mathematics, University of Tabriz, Tabriz, Iran;
and School of Mathematics, Institute for Research in Fundamental
Sciences (IPM), P.O. Box. 19395-5746, Tehran, Iran.}
\email{naghipour@ipm.ir} \email {naghipour@tabrizu.ac.ir}
\address{Department of Mathematics, Azarbaijan Shahid Madani University, Tabriz, Iran. }%
\email{m\_sedghi@tabrizu.ac.ir} \email {sedghi@azaruniv.ac.ir}%
\thanks{ 
$^*$The second author is grateful to the hospitality and facilities offered of the Max-Planck Institut f\"ur Mathematik (Bonn) during the preparation of this paper.\\
$^\dag$Corresponding author: e-mail: {\it naghipour@ipm.ir} (Reza Naghipour)}%
%\subjclass{}%
%\date{}%
%\dedicatory{}%
%\commby{}%
% ----------------------------------------------------------------
\begin{abstract}
The purpose of the present paper is to continue the study of modules cofinite and weakly cofinite  with respect to an ideal $\frak a$ of a Noetherian ring $R$. It is shown that an $R$-module $M$ is cofinite  with respect to  $\frak a$, if and only if, $\Ext^i_R(R/\frak a,M)$ is finitely generated for all $i\leq {\rm cd}(\frak a,M)+1$, whenever $\dim R/\frak a=1$.
In addition, we show that if $M$ is finitely generated and  $H^i_{\frak a}(M)$
are weakly Laskerian for all $i\leq t-1$, then $H^i_{\frak a}(M)$ are ${\frak a}$-cofinite for all $i\leq t-1$ and for any minimax submodule
$K$  of  $H^{t}_{\frak a}(M)$,  the $R$-modules $\Hom_R(R/{\frak a}, H^{t}_{\frak a}(M)/K)$ and $\Ext^{1}_R(R/{\frak a}, H^{t}_{\frak a}(M)/K)$
are finitely generated, where $t$ is a non-negative integer.  Finally, we explore a criterion for weakly cofiniteness of modules  with respect to an ideal of dimension one. Namely for such ideals it suffices that the two first $\Ext$-modules in the definition for weakly cofiniteness are weakly Laskerian. As an application of this result we deduce  that the category of all ${\frak a}$-weakly cofinite modules over  $R$ forms a full Abelian subcategory of the category of modules.

\end{abstract}
\maketitle
% ----------------------------------------------------------------
\section{Introduction}
Let  $R$  denote a commutative Noetherian ring
(with non-zero identity) and $\frak a$ an ideal of $R$.  Also, we let $M$  denote an arbitrary $R$-module.
  
It is well-known result that if $R$ is a local (Noetherian) ring
with maximal ideal $\frak m$, then the $R$-module $M$ is Artinian
if and only if $\Supp(M)\subseteq \{\frak m\}$ and $\Ext^{j}_{R}(R/\frak m,M)$ is finitely generated for
all $j\geq0$ (cf. \cite[Proposition 1.1]{Ha}).

Using this idea, Hartshorne \cite{Ha} introduced the class of cofinite modules, answering in negative
a question of Grothendieck (cf. \cite[Expos$\acute{e}$ XIII, Conjecture 1.1]{Gr2}).
In fact, Grothendieck conjectured that for any ideal $\frak a$ of $R$
and any finitely generated $R$-module $M$, the $R$-module $\Hom_R(R/\frak a,H^{i}_{\frak a}(M))$ is finitely generated,
where $H^{i}_{\frak a}(M)$ is the $i$-th local cohomology module of $M$ with support in $V(\frak a)$, (this is the case
when $\frak a=\frak m$, the maximal ideal in a local ring, since the modules $H^{i}_{\frak m}(M)$ are Artinian), but soon
Hartshorne was able to present a counterexample (see \cite{Ha} for details and proof) which shows that
this conjecture is false even when $R$ is regular, and where he defined an $R$-module $M$ to be {\it cofinite with respect to} $\frak a$ (abbreviated as  $\frak a$-{\it cofinite})
if the support of  $M$ is contained in $V(\frak a)$ and $\Ext^{j}_{R}(R/\frak a,M)$ is
finitely generated for all $j$ and asked the following questions:

(i) {\it For which rings $R$ and ideals $\frak a$ are the modules
$H^{i}_{\frak a}(M)$,  $\frak a$-cofinite for all $i$ and all finitely generated
modules $M$?}\\
(ii) {\it Whether the category $\mathscr{C}(R, \frak a)_{cof}$
of $\frak a$-cofinite modules forms an Abelian subcategory of the category of all $R$-modules?}\\

With respect to the question (i), Hartshorne in \cite{Ha} and later
Chiriacescu in \cite{Ch} showed that if $R$ is a complete regular
local ring and $\frak a$ is a prime ideal such that $\dim R/\frak a=1$, then
$H^{i}_{\frak a}(M)$ is $\frak a$-cofinite for any finitely generated $R$-module
$M$ (see \cite[Corollary 7.7]{Ha}).

Also, Delfino and Marley \cite[Theorem 1]{DM} and Yoshida \cite
[Theorem 1.1]{Yo} have eliminated the complete hypothesis entirely.
Finally, more recently Bahmanpour and Naghipour removed the local
condition on the ring (see \cite[Theorem 2.6]{BN1}).

 For a survey of recent developments on finiteness properties of local cohomology modules, see Lyubeznik's interesting paper \cite{Ly3}.

In the second section, we establish several characterizations of the $\frak a$-cofiniteness of an $R$-module $M$. More precisely we prove the following result:

\begin{thm}
Let $R$ be a Noetherian ring, $M$ an $R$-module and $\frak a$ a one-dimensional ideal of $R$ such that
 $\Supp(M)\subseteq V(\frak a)$. Then the following conditions are equivalent:

{\rm(i)} $M$ is $\frak a$-cofinite. 

{\rm(ii)} $H^i_{\frak a}(M)$ is $\frak a$-cofinite, for all $i$.

{\rm(iii)} $\Ext^i_R(R/\frak a,M)$ is finitely generated, for all $i\leq {\rm cd}(\frak a,M)+1.$

{\rm(iv)} $\Ext^i_R(N,M)$ is finitely generated,  for all $i\leq {\rm cd}(\frak a,M)+1$ and for any finitely generated
$R$-module $N$ with $\Supp(N) \subseteq V(\frak a)$.

{\rm(v)} $\Ext^i_R(N,M)$ is finitely generated, for all $i\leq {\rm cd}(\frak a,M)+1$ and for some finitely generated $R$-module $N$ with $\Supp(N) = V(\frak a)$.
\end{thm}

Pursuing this point of view further we derive the following consequence of Theorem 1.1, which is an extension of the main results of Delfino-Marley \cite{DM} and Yoshida \cite{Yo} for an arbitrary Noetherian ring $R$.

\begin{cor}
Let $R$ be a Noetherian ring and let $\frak a, \frak b$ be ideals of $R$ such that
$\frak b\subseteq \Rad(\frak a)$. Let $M$ be a $\frak b$-cofinite  $R$-module.

${\rm (i)}$ If $\dim R/\frak a=1$, then  $H^i_{\frak a }(M)$ is $\frak a$-cofinite  for all $i$.

${\rm (ii)}$ If $\dim R/\frak b=1$, then  $H^i_{\frak b}(M)$ is $\frak a$-cofinite for all $i$.
\end{cor}
In \cite{Zo1} H. Z\"{o}schinger, introduced the interesting class of
minimax modules, and he has in \cite{Zo1, Zo2} given many equivalent
conditions for a module to be minimax. The $R$-module $N$ is said to
be {\it minimax}, if there is a finitely generated
submodule $L$ of $N$, such that $N/L$ is Artinian. The class of
minimax modules thus includes all finitely generated and all
Artinian modules. It was shown by T. Zink \cite{Zi} and by E. Enochs
\cite{En} that a module over a complete local
ring is minimax if and only if it is Matlis reflexive.
 
In the second section, we also shall prove the following, which is a generalization of the main result of Brodmann-Lashgari \cite{BL}.
\begin{thm}
Let $R$ be a  Noetherian ring, $\frak a$ an ideal of $R$ and $M$ a
finitely generated $R$-module such that for a non-negative integer $t$, the $R$-modules $H^i_{\frak a}(M)$
are weakly Laskerian for all $i\leq t$. Then the $R$-modules
$H^0_{\frak a}(M),\dots,H^t_{\frak a}(M)$ are ${\frak a}$-cofinite and for any minimax submodule
$K$  of  $H^{t+1}_{\frak a}(M)$ and for any finitely generated
$R$-module $L$ with $\Supp(L)\subseteq V({\frak a})$,  the $R$-modules $\Hom_R(L,H^{t+1}_{\frak a}(M)/K)$ and $\Ext^{1}_R(L,H^{t+1}_{\frak a}(M)/K)$
are finitely generated.
\end{thm}

An $R$-module $M$ is said to be a {\it weakly Laskerian module}, if the set of associated primes of any quotient module $M$ is finite (see \cite{DM1}  and \cite{Ro}).

With respect to the question (ii), Hartshorne with an example showed that this not true in general. However, he proved that if $\frak a$ is a prime ideal of dimension one in a complete regular local ring $R$, then the answer to his question
is yes. In \cite{DM}, Delfino and Marley extended this result to arbitrary
complete local rings. Recently, Kawasaki \cite{Ka2}, by using a spectral sequence argument, generalized the Delfino and Marley's result
for an arbitrary ideal $\frak a$ of dimension one in a local ring $R$. Finally, more recently Bahmanpour, Naghipour and Sedghi in \cite{BNS1} removed the local
condition on the ring. Namely, therein it is shown that Hartshorne's question is true
for $\mathscr{C}^1(R, \frak a)_{cof}$,  the category of all $\frak a$-cofinite
$R$-modules $M$ with $\dim\Supp(M)\leq 1$,  for all ideals $\frak a$ in a  Noetherian ring $R$. The proof of this result
is based on \cite[Proposition 2.6]{BNS1} which states that in order to deduce the $\frak a$-cofiniteness  for a module $M$ with
$\dim\Supp(M)\leq 1$ and $\Supp(M)\subseteq V(\frak a)$, it suffices that we know that the $R$-modules $\Hom_R(R/\frak a,M)$ and $\Ext^1_R(R/\frak a,M)$
are finitely generated.

The main goal of Section 3 is to establish  the analogue of this result to the $\frak a$-weakly cofiniteness. Namely, in this section among other things,  we show that for the $\frak a$-weakly cofiniteness  of a module $M$ with
$\dim\Supp(M)\leq 1$ and $\Supp(M)\subseteq V(\frak a)$, it suffices that we know that the $R$-modules $\Hom_R(R/\frak a,M)$ and $\Ext^1_R(R/\frak a,M)$ are weakly Laskerian. In particular, when $\frak a$
is one-dimensional, in order to deduce the $\frak a$-weakly cofiniteness for a module (with support in $V(\frak a)$), it suffices that we know that the first two $\Ext$-modules in the definition for weakly cofiniteness are weakly Laskerian. More precisely, we shall show that:

\begin{thm}
Let $\frak a$ denote an ideal of a Noetherian ring $R$ and let $M$ be an
$R$-module such that $\dim\Supp(M)\leq 1$ and $\Supp(M) \subseteq V(\frak a)$. Then $M$ is $\frak a$-weakly cofinite if and only if
the $R$-modules ${\rm Hom}_R(R/\frak a,M)$ and ${\rm Ext}^1_R(R/\frak a,M)$
are weakly Laskerian.
\end{thm}
 An $R$-module  $M$ is said to be $\frak a$-{\it weakly cofinite} if
$\Supp(M)\subseteq V(\frak a)$ and ${\rm Ext}^{i}_{R}(R/\frak a, M)$ is a weakly Laskerian module for all $i$ (see  \cite{DM2}). We denote the category of the $\frak a$-weakly cofinite modules by $\mathscr{C}(R, \frak a)_{wcof}$. As an application of Theorem 1.4  we show that,  when $\frak a$ is one-dimensional,   $\mathscr{C}(R, \frak a)_{wcof}$ forms an Abelian subcategory of the category of all $R$-modules (see Corollary 3.6). That is, if $f: M\longrightarrow N$ is an $R$-homomorphism between  $\frak a$-weakly cofinite modules, then $\ker f$ and ${\rm coker} f$ are $\frak a$-weakly cofinite.  The proof of this result is based on the following theorem. 
\begin{thm}
Let $\frak a$ be an ideal of a Noetherian ring $R$. Let
$\mathscr{C}^1(R, \frak a)_{wcof}$ denote the category of $\frak a$-weakly cofinite
$R$-modules $M$ with $\dim\Supp(M)\leq 1$. Then
$\mathscr{C}^1(R, \frak a)_{wcof}$ is an Abelian category.
\end{thm}

The proof of Theorem 1.5 is given in Theorem 3.5.  Finally, we end the paper with a question concerning the Serre subcategory.\\

Throughout this paper, $R$ will always be a commutative Noetherian
ring with non-zero identity and $\frak a$ will be an ideal of $R$. For an $R$-module $M$, the
$i$-th local cohomology module of $M$ with support in $\frak a$ is
defined as
 $$H^i_{\frak a}(M) = \underset{n\geq1} {\varinjlim}\,\, \text{Ext}^i_R(R/\frak a^n, M).$$ 
 For facts about the local cohomology modules we refer to the textbook by Brodmann-Sharp \cite{BS} or Grothendieck's interesting book \cite{Gr1}.

Further, for any ideal $\frak b$ of $R$, we denote the set
$\{\frak p \in {\rm Spec}\,R:\, \frak p\supseteq \frak b \}$ by
$V(\frak b)$; and  the {\it radical} of $\frak{b}$, denoted by $\Rad(\frak{b})$, we define to be the set $\{x\in R \,: \, x^n \in \frak{b}$ for some $n \in\mathbb{N}\}$.

 For an Artinian $R$-module $A$ the set of attached
prime ideals of $A$ is denoted by $\Att_R A$. Also, for each $R$-module $L$, we denote by
 ${\rm Ass h}_RL$  the set $\{\frak p\in \Ass
_RL:\, \dim R/\frak p= \dim L\}$. Finally, we shall use $\Max(R)$ to denote the set of all maximal
ideals of $R$.  For any unexplained notation and terminology we refer
the reader to \cite{BH} and \cite{Mat}. \\

%In \cite{Zo1} H. Z\"{o}schinger, introduced the interesting class of minimax modules, and he has in \cite{Zo1, Zo2} given many equivalent
%conditions for a module to be minimax. The $R$-module $N$ is said to be a {\it minimax module}, if there is a finitely generated
%submodule $L$ of $N$, such that $N/L$ is Artinian. The class of minimax modules thus includes all finitely generated and all
%Artinian modules. It was shown by T. Zink \cite{Zi} and by E. Enochs \cite{En} that a module over a complete local ring is minimax if and only if it is Matlis %reflexive. \\

%\end{align*}
\section{ Modules cofinite}
The main goals of this section are Theorems 2.4 and 2.8.  The following lemmas will be needed in the proof of these results. Recall that a class $\mathcal{S}$ of $R$-modules is a {\it Serre subcategory} of the category of $R$-modules, when it is closed under taking submodules, quotients and extensions.
It is well known  that the subcategories of, finitely generated, minimax, weakly Laskerian, and
Matlis reflexive modules are examples of Serre subcategory. Following we let $\mathcal{S}$ denote
a Serre subcategory of the category of $R$-modules.

\begin{lem}
    \label{2.1}
 Let $R$ be a Noetherian ring and $\frak a$ an ideal of $R$. Let $s$  be a
 non-negative integer and let $M$ be an $R$-module such that $\Ext_R^s(R/{\frak a},M)\in \mathcal{S}$.
Suppose that $\Ext_R^j(R/\frak a,H^i_{\frak a}(M))\in \mathcal{S}$ for all $i<s$ and all $j\geq0$. Then
$\Hom_R(R/\frak a,H^s_{\frak a}(M))\in \mathcal{S}$.
\end{lem}
\proof See \cite[Theorem 2.2]{AS}.\qed

\begin{lem}
    \label{2.2}
Let $R$ be a Noetherian ring and ${\frak a}$ an ideal of $R$. Let $s$  be a
 non-negative integer and let $M$ be an $R$-module such that $\Ext_R^{s+1}(R/{\frak a},M)\in \mathcal{S}$.
Suppose that $\Ext_R^j(R/{\frak a},H^i_{\frak a}(M))\in \mathcal{S}$ for all $i<s$ and all $j\geq0$. Then
$\Ext^1_R(R/{\frak a},H^s_{\frak a}(M))\in \mathcal{S}$.
\end{lem}
\proof We use induction on $s$. Let $s=0$. Then the exact sequence
$$0\longrightarrow \Gamma_{\frak a}(M)\longrightarrow M\longrightarrow M/\Gamma_{\frak a}(M)\longrightarrow0,\,\,\,\,\,\,\,\,\,\,\,\,\,\, (\dag)$$
induces the exact sequence
$$\Hom_R(R/{\frak a},M/\Gamma_{\frak a}(M))\longrightarrow \Ext^1_R(R/{\frak a}, \Gamma_{\frak a}(M)) \longrightarrow \Ext^1_R(R/{\frak a}, M).$$
 As $\Hom_R(R/\frak a,M/\Gamma_{\frak a}(M))=0$ and $\Ext_R^{1}(R/{\frak a},M)$ are in $\mathcal{S}$, it follows that
$$\Ext_R^{1}(R/{\frak a},\Gamma_{\frak a}(M))\in \mathcal{S}.$$

Now, suppose inductively that $s>0$ and that the assertion holds for $s-1$. Using the exact sequence $(\dag)$
 we obtain the following exact sequence, $j\geq0$,
$$\Ext^j_R(R/{\frak a},M)\longrightarrow \Ext^j_R(R/{\frak a}, M/\Gamma_{\frak a}(M)) \longrightarrow \Ext^{j+1}_R(R/{\frak a}, \Gamma_{\frak a}(M)).$$
Therefore, since $\Ext^{s+2}_R(R/{\frak a}, \Gamma_{\frak a}(M))$ and $\Ext^{s+1}_R(R/{\frak a},M)$ are in $\mathcal{S}$, it follows that
 $\Ext^{s+1}_R(R/{\frak a}, M/\Gamma_{\frak a}(M))\in\mathcal{S}$. Also, it  easily follows from assumption and \cite[Corollary 2.1.7]{BS}
that $\Ext_R^j(R/{\frak a},H^i_{\frak a}(M/\Gamma_{\frak a}(M)))\in \mathcal{S}$ for all $i<s$ and all $j\geq0$. Therefore we may
assume that $\Gamma_{\frak a}(M)=0$.

Next, let $E_R(M)$ denote the injective hull of $M$. Then $\Gamma_{\frak a}(E_R(M))=0$,
and so it follows from the exact sequence
$$0\longrightarrow M\longrightarrow E_R(M)\longrightarrow E_R(M)/M\longrightarrow 0,$$
that $H^{i+1}_{\frak a}(M)\cong H^{i}_{\frak a}(E_R(M)/M)$ for all $i\geq0$. Also, as $\Hom_R(R/{\frak a}, E_R(M))=0$,
it yields that $$\Ext_R^j(R/{\frak a},M)\cong \Ext_R^{j+1}(R/{\frak a}, M),$$ for all $j\geq0$.
Consequently the $R$-module $E_R(M)/M$ satisfies our condition hypothesis. Thus $\Ext^1_R(R/{\frak a},H^{s-1}_{\frak a}(E_R(M)/M))\in \mathcal{S}$.
Now the assertion follows from the isomorphism $$H^{s}_{\frak a}(M)\cong H^{s-1}_{\frak a}(E_R(M)/M).$$ \qed

\begin{lem}
 \label{2.3}
 Let ${\frak a}$ be an ideal of a Noetherian ring $R$ and $M$  a non-zero
$R$-module, such that $\dim\Supp(M)\leq 1$ and $\Supp(M)\subseteq V({\frak a})$. Then the
following statements are equivalent:

{\rm(i)} $M$ is ${\frak a}$-cofinite.

{\rm(ii)} The $R$-modules $\Hom_R(R/{\frak a},M)$ and $\Ext^1_R(R/{\frak a},M)$ are finitely generated.
\end{lem}
\proof  See \cite[Proposition 2.6]{BNS1}. \qed\\

Now we are prepared to state and prove  the first main theorem of this section. Recall that for an $R$-module $N$, the {\em cohomological
 dimension of } $N$ {\em with respect to} an ideal $\frak a$ of $R$,
 denoted by $\text{cd}(\frak a, N)$, is defined as
 $$\text{cd}(\frak a, N)= \text{sup}\{i\in \mathbb N_0\mid H^i_{\frak a}(N) \neq 0\}.$$

\begin{thm}
 \label{2.4}
Let $R$ be a Noetherian ring, $M$ an  $R$-module and ${\frak a}$ a one-dimensional ideal of $R$. Then the following conditions are equivalent:

{\rm(i)} $\Ext^i_R(R/{\frak a},M)$ is finitely generated, for all $i\leq {\rm cd}({\frak a},M)+1.$

{\rm(ii)} $H^i_{\frak a}(M)$ is ${\frak a}$-cofinite,  for all $i$.

{\rm(iii)} $\Ext^i_R(R/{\frak a},M)$ is finitely generated,  for all $i$.

{\rm(iv)} $\Ext^i_R(N,M)$ is finitely generated, for all $i\leq {\rm cd}({\frak a},M)+1$ and for any finitely generated
$R$-module $N$ with $\Supp(N) \subseteq V({\frak a})$.

{\rm(v)} $\Ext^i_R(N,M)$ is finitely generated, for all $i\leq {\rm cd}({\frak a},M)+1$ and for some finitely generated
$R$-module $N$ with $\Supp(N) = V({\frak a})$.

{\rm(vi)} $\Ext^i_R(N,M)$ is finitely generated, for all $i$ and for any finitely generated
$R$-module $N$ with $\Supp(N) \subseteq V({\frak a})$.

{\rm(vii)} $\Ext^i_R(N,M)$ is finitely generated, for all $i$ and for some finitely generated
$R$-module $N$ with $\Supp(N) = V({\frak a})$.
\end{thm}

\proof In order to prove ${\rm(i)}\Longrightarrow{\rm(ii)}$ we may assume that $i\leq {\rm cd}({\frak a},M)$. Now, we use induction on $i$.   When $i=0$, then the exact sequence
$$0\longrightarrow \Gamma_{\frak a}(M)\longrightarrow M\longrightarrow M/\Gamma_{\frak a}(M)\longrightarrow0,$$
induces the exact sequence
$$0\longrightarrow \Hom_R(R/{\frak a},\Gamma_{\frak a}(M))\longrightarrow\Hom_R(R/{\frak a}, M)\longrightarrow \Hom_R(R/{\frak a},M/\Gamma_{\frak a}(M))$$$$\longrightarrow \Ext^1_R(R/{\frak a}, \Gamma_{\frak a}(M)) \longrightarrow \Ext^1_R(R/{\frak a}, M).$$

 As $\Hom_R(R/{\frak a},M/\Gamma_{\frak a}(M))=0$ and $\Ext_R^{j}(R/{\frak a},M)$, for $j=0, 1$,  is finitely generated, it follows that
$\Hom_R(R/{\frak a}, \Gamma_{\frak a}(M)))$ and $\Ext_R^{1}(R/{\frak a},\Gamma_{\frak a}(M))$ are finitely generated. It now follows from Lemma 2.3 that
 $\Gamma_{\frak a}(M)$ is ${\frak a}$-cofinite.

 Assume, inductively, that $i>0$ and that the result has been proved for $i-1$. Then the $R$-modules  $$H^0_{\frak a}(M),H^1_{\frak a}(M), \dots, H^{i-1}_{\frak a}(M),$$
are ${\frak a}$-cofinite, and so it follows from Lemmas 2.1 and 2.2 that  $\Hom_R(R/{\frak a},H^i_{\frak a}(M))$ and $\Ext^1_R(R/{\frak a},H^i_{\frak a}(M))$
are finitely generated. Now, it yields from  Lemma 2.3  that
 $H^i_{\frak a}(M)$ is ${\frak a}$-cofinite.

  The implication ${\rm(ii)}\Longrightarrow{\rm(iii)}$ follows from \cite[Proposition 3.9]{Me2}, and for prove ${\rm(iii)}\Longrightarrow{\rm(vi)}$  see \cite[Lemma 1]{Ka1}. Finally, in order to
complete the proof, it is enough for us to show that  ${\rm(v)}\Longrightarrow{\rm(iv)}$. To this end, let $L$ be a finitely generated
$R$-module with $\Supp(L) \subseteq V({\frak a})$ and $N$ a  finitely generated $R$-module  such that $\Supp(N) = V({\frak a})$.
Then $\Supp(L)\subseteq \Supp(N)$, and so according to
Gruson's Theorem
 \cite[Theorem 4.1]{Va}, there exists a chain \[0=L_0\subset L_1\subset\cdots\subset L_k=L,\]
 such that the factors $L_j/L_{j-1}$ are homomorphic images of a direct sum of finitely
 many copies of $N$. Now consider the exact sequences
\[0\longrightarrow K\longrightarrow
  N^n\longrightarrow L_1\longrightarrow 0  \ \ \ \ \ \ \ \ \ \ \ \ \ \ \ \ \ \ \      \]

\[0\longrightarrow L_1\longrightarrow
   L_2\longrightarrow L_2/L_1\longrightarrow 0  \ \ \ \ \ \ \ \ \ \ \ \ \ \ \ \ \ \ \      \]
                               $$\vdots$$
                                       \[0\longrightarrow L_{k-1}\longrightarrow
  L_k\longrightarrow L_k/L_{k-1}\longrightarrow 0,  \ \ \ \ \ \ \ \ \ \ \ \ \ \ \ \ \ \ \  \]

for some positive integer $n$.  Now,  from the long exact sequence

 \[\cdots\rightarrow {\rm Ext}_R^{i-1}(L_{j-1},N)\rightarrow
   {\rm Ext}_R^{i}(L_j/L_{j-1},N)\rightarrow {\rm
   Ext}_R^{i}(L_j,N) \rightarrow {\rm Ext}_R^i(L_{j-1},N)\rightarrow\cdots,\]

   and an easy induction on $k$, it suffices  to prove the case when
   $k=1$.

 Thus there is an exact sequence \[0\longrightarrow K\longrightarrow
  N^n\longrightarrow L\longrightarrow 0  \ \ \ \ \ \ \ \ \ \ \ \ \ \ \ \ \ \ \ \ \ \ \ \ \ \ \ \ \,\,\,\,\,\,\,\,      (\dag)\]   for some $n\in\Bbb N$
and some finitely generated $R$-module $K$.

    Now, we use induction on $i$. First, $\Hom_R(L,M)$ is a submodule
     of $\Hom_R(N^n,M)$; hence, in view of assumption,  $\Ext_R^0(L,M)$ is finitely
  generated. So assume that $i>0$ and that $\Ext_R^j(L^\prime, M)$
     is finitely   generated for every finitely
  generated $R$-module $L^\prime$ with $\Supp(L')\subseteq \Supp(N)$ and for
  all $j\leq i-1$. Now, the exact sequence $(\dag)$ induces the long
  exact sequence \[\cdots\longrightarrow \Ext_R^{i-1}(K,M)\longrightarrow
   \Ext_R^{i}(L,M)\longrightarrow \Ext_R^{i}(N^n,M)\longrightarrow\cdots,\] so that, by the inductive
   hypothesis, $\Ext_R^{i-1}(K,M)$ is finitely generated. On the
   other hand $\Ext_R^{i}(N^n,M)\cong \stackrel{n}{\oplus}\Ext_R^{i}(L,M)$
   is  finitely generated, and so $\Ext_R^{i}(L,M)$ is finitely
  generated,  the inductive step is complete.  \qed\\

As a consequence of Theorem 2.4, we derive the following result which is an extension of the main results of Delfino-Marley \cite{DM} and Yoshida \cite{Yo} for arbitrary Noetherian rings.
\begin{cor}
Let $R$ be a Noetherian ring and ${\frak a}, {\frak b}$ be ideals of $R$ such that
${\frak b}\subseteq \Rad({\frak a})$. Let $M$ be a ${\frak b}$-cofinite $R$-module.

${\rm (i)}$ If $\dim R/{\frak a}=1$, then the $R$-module $H^i_{\frak a}(M)$ is ${\frak a}$-cofinite for all $i$.

${\rm (ii)}$ If $\dim R/{\frak b}=1$, then the $R$-module $H^i_{\frak b}(M)$ is ${\frak a}$-cofinite for all $i$.
\end{cor}

\proof In order to show (i), since ${\frak b}\subseteq \Rad({\frak a})$, it follows that $\Supp(R/{\frak a})\subseteq \Supp (R/{\frak b})$. On the other hand,
since $M$ is ${\frak b}$-cofinite it follows from \cite[Lemma 1]{Ka1}
that $M$ is also ${\frak a}$-cofinite.  Now as $\dim R/{\frak a}=1$,  it follows from Theorem 2.4 that  $H^i_{\frak a}(M)$ is ${\frak a}$-cofinite for all $i$.

To prove (ii), since $\dim R/{\frak b}=1$ and $M$ is ${\frak b}$-cofinite it follows from Theorem 2.4 that
$H^i_{\frak b}(M)$ is ${\frak b}$-cofinite for all $i$. Now, because of $\Supp(R/{\frak a})\subseteq \Supp (R/{\frak b})$ it
follows from \cite[Lemma 1]{Ka1}  that $H^i_{\frak b}(M)$ is ${\frak a}$-cofinite, for all $i$. \qed\\

Before proving  the next main theorem, we need the following lemma and proposition, which will be used in Theorem 2.8.

\begin{lem}
    \label{2.7}
Let $R$ be a Noetherian ring and $M$ an $R$-module. Then $M$ is weakly Laskerian
if and only if there exists a finitely generated submodule $N$ of $M$ such that
$\Supp(M)/N$ is finite.
\end{lem}
\proof See \cite[Theorem 3.3]{Ba1}.\qed\\

  \begin{prop}
 \label{2.8}
Let $R$ be a Noetherian ring, ${\frak a}$ an ideal of $R$ and $M$ a
finitely generated $R$-module such that  $H^i_{\frak a}(M)$
is weakly Laskerian for all $i\leq t$. Then the $R$-modules
$$H^0_{\frak a}(M),\dots,H^t_{\frak a}(M)$$ are ${\frak a}$-cofinite. In addition the $R$-modules 
\begin{center}
$\Hom_R(R/{\frak a},H^{t+1}_{\frak a}(M))$ and $\Ext^{1}_R(R/{\frak a},H^{t+1}_{\frak a}(M))$
\end{center}
 are finitely generated. In particular, the  set $\Ass_R H_{\frak a}^{t+1}(M)$ is finite.
\end{prop}
\proof We use induction on $t$. The case $t=0$ follows from
Lemmas 2.1 and 2.2. So, let $t\geq 1$ and the case $t-1$ is settled. Then by
inductive hypothesis the $R$-modules $H^0_{\frak a}(M),\dots,H^{t-1}_{\frak a}(M)$
are ${\frak a}$-cofinite and the $R$-modules 
\begin{center}
$\Hom_R(R/{\frak a},H^{t}_{\frak a}(M))$ and 
$\Ext^{1}_R(R/{\frak a},H^{t}_{\frak a}(M))$ 
\end{center}
are finitely generated. Now
since by assumption the $R$-module $H^t_{\frak a}(M)$ is weakly Laskerian,
it follows from Lemma 2.6 that there is a finitely generated
submodule $N$ of $H^t_{\frak a}(M)$ such that $\Supp(H^t_{\frak a}(M)/N)$ is finite
set, and so $\dim\Supp (H^t_{\frak a}(M)/N)\leq 1$. Now it follows from the
exact sequence $$0\longrightarrow N \longrightarrow H^t_{\frak a}(M)
\longrightarrow H^t_{\frak a}(M)/N \longrightarrow 0,$$  that the
$R$-modules $$\Hom_R(R/{\frak a},H^t_{\frak a}(M)/N)\,\,\,\,\,{\rm
and}\,\,\,\,\,\Ext^1_R(R/{\frak a},H^t_{\frak a}(M)/N),$$ are finitely generated.
Therefore it follows from Lemma 2.3 that the $R$-module
$H^t_{\frak a}(M)/N$ is ${\frak a}$-cofinite, and so the $R$-module
$H^t_{\frak a}(M)$ is ${\frak a}$-cofinite. Hence, it follows from Lemmas 2.1 and 2.2 that the
$R$-modules ${\rm Hom}_R(R/{\frak a},H^{t+1}_{\frak a}(M))$ and ${\rm
Ext}^{1}_R(R/{\frak a},H^{t+1}_{\frak a}(M))$ are finitely generated. This completes
the induction step. \qed\\

Now, we are ready to state and prove the second main result of this section, which is a generalization the main results of Bahmanpour-Naghipour \cite[Theorem 2.6]{BN} and Brodmann-Lashgari \cite[Theorem 2.2]{BL}.
\begin{thm}
Let $R$ be a  Noetherian ring, $\frak a$ an ideal of $R$ and $M$ a
finitely generated $R$-module such that for a non-negative integer $t$, the $R$-modules $H^i_I(M)$
are weakly Laskerian for all $i\leq t$. Then the $R$-modules
$$H^0_{\frak a}(M),\dots,H^t_{\frak a}(M)$$ are $\frak a$-cofinite and for any minimax submodule
$K$ of $H^{t+1}_{\frak a}(M)$ and for any finitely generated
$R$-module $L$ with $\Supp(L)\subseteq V({\frak a})$,  the $R$-modules
\begin{center}
 $\Hom_R(L, H^{t+1}_{\frak a}(M)/K)$ and $\Ext^{1}_R(L, H^{t+1}_{\frak a}(M)/K)$
\end{center}
are finitely generated.
\end{thm}
\proof By virtue of Proposition 2.7 the $R$-module $H^i_{\frak a}(M)$ is $\frak a$-cofinite for all $i\leq t$ and 
$\Hom_R(R/{\frak a},H^{t+1}_{\frak a}(M))$ is finitely generated. Hence the $R$-module
$\Hom_R(R/{\frak a}, K)$ is  finitely generated, and so in view of \cite[Proposition 4.3]{Me2}, $K$ is ${\frak a}$-cofinite.
Thus, \cite[Lemma 1]{Ka1} implies that  $\Ext_R^{i}(L,K)$ is  finitely generated for all $i$.

Next, the exact sequence
\[0\longrightarrow {K}\longrightarrow
   H_{\frak a}^{t+1}(M)\longrightarrow H_{\frak a}^{t+1}(M)/K \longrightarrow 0\]
provides the following exact sequence,

\[{\Hom_R(L,{\rm H}_{\frak a}^{t+1}(M))}\longrightarrow
   \Hom_R(L, H_{\frak a}^{t+1}(M)/K)\longrightarrow  \Ext_R^{1}(L,K)$$ $$\longrightarrow \Ext^1_R(L, H_{\frak a}^{t+1}(M))\longrightarrow \Ext^1_R(L, H_{\frak a}^{t+1}(M)/K)\longrightarrow    \Ext_R^{2}(L,K) .\]

Now,  since  $\Ext_R^{i}(L, K)$ is finitely generated, the assertion follows from  Proposition 2.7 and \cite[Lemma 1]{Ka1}, because the $R$-modules  
\begin{center}
$\Hom_R(L, H_{\frak a}^{t+1}(M))$ and $\Ext^1_R(L, H_{\frak a}^{t+1}(M))$
\end{center}
are finitely generated.   \qed\\

\section{ Modules weakly cofinite}
The purpose of this section is to establish that the category of modules weakly cofinite with respect to an ideal of dimension one in a Noetherian ring is a full Abelian
subcategory of the category of modules. The main goal of this section is Theorem 3.5. The proof of this theorem  is based on the Proposition 3.2, which plays a key role in this section, says that (when $\frak a$ is one-dimensional), in order to deduce the $\frak a$-weakly cofiniteness for a module (with support in $V(\frak a)$), it suffices that we know that the first two $\Ext$-modules in the definition for weakly cofiniteness are weakly Laskerian. Before stating it, we record a  lemma that will be needed in the proof of this proposition.

%\begin{lem}
   % \label{2.1}
 %Let $R$ be a Noetherian ring and $\frak a$ an ideal of $R$.
%Then, for any $R$-module $T$, the following conditions are equivalent:

%${\rm(i)}$ $\Ext_R^n(R/\frak a,T)$ is weakly Laskerian for all $n\geq0$,

%${\rm(ii)}$ $\Ext_R^n(N,T)$ is  weakly Laskerian for all $n\geq0$ and for each
%finitely generated $R$-module $N$ for which $\Supp(N)\subseteq V(\frak a)$.
%\end{lem}
%\proof See \cite[Lemma 2.8]{DM2}.\qed\\

%\begin{lem}
   % \label{2.2}
%Let $R$ be a Noetherian ring and $M$ an $R$-module. Then $M$ is weakly Laskerian
%if and only if there exists a finitely generated submodule $N$ of $M$ such that $\Supp(M)/N$ is finite.
%\end{lem}
%\proof See \cite[Theorem 3.3 ]{Ba1}.\qed

%\begin{lem}
    %\label{2.3}
%Let $R$ be a Noetherian ring and let $M, N$ be two
%$R$-modules such that $M$ is weakly Laskerian and $N$ is finitely generated.
%Then $\Ext_R^i(N, M)$ is weakly Laskerian for all $i\geq0$.
%\end{lem}
%\proof See \cite[Lemma 2.2 ]{DM2}.\qed\\

\begin{lem}
 \label{sam 2.4}
Let $(R,\frak m)$  be a local (Noetherian) ring and let $A$ be an Artinian $R$-module.\\
$\rm(i)$ If $\frak a$ is an ideal of $R$ such that   $\Hom_{R}(R/\frak a,A)$ is a
finitely generated $R$-module, then $$V(\frak a)\cap \Att_{R} A\subseteq V(\mathfrak{m}).$$
$\rm(ii)$ If  $x$ is an element of $R$ such that $V(Rx)\cap\Att_{R}A\subseteq\{\mathfrak{m}\}$, then the $R$-module $A/xA$ has finite length.
\end{lem}
\proof See \cite[Lemmas 2.4 and 2.5]{BN1}.\qed\\

%\begin{lem}
%Let $(R,\mathfrak{m})$ and $A$ be as in Lemma {\rm 2.4.} Suppose that $x$ is an
%element in $\mathfrak{m}$ such that $V(Rx)\cap \Att_{R}A\subseteq\{\mathfrak{m}\}$. Then the $R$-module $A/xA$ has finite length.
%\end{lem}
%\proof See \cite[Lemma 2.4 ]{BN1}.\qed\\
%\begin{lem}
 %   \label{2.3}
%Let $R$ be a Noetherian ring and $I$ an ideal of $R$. Let $M$ be a minimax
%$R$-module such that $\Supp(M)\subseteq V(I)$. Then $M$ is $I$-cofinite
 % if and only if $(0:_M I)$ is finitely generated.
%\end{lem}
%\proof See \cite[Proposition 4.3 ]{Me2}.\qed\\

The following proposition will be one our main tools in this section. It's proof  is based on the important notion of the arithmetic rank
of an ideal. The {\it arithmetic rank} of an ideal ${\frak b}$ in a
 Noetherian ring $R$, denoted by ${\rm ara}({\frak b})$, is the
least number of elements of $R$ required to generate an ideal which has the
same radical as ${\frak b}$, i.e.,
$${\rm ara}({\frak b})=\min\{n\in \mathbb{N}_0: \exists b_1, \dots, b_n \in R\,\,  \text{with}\,\, \Rad(b_1, \dots, b_n)=\Rad({\frak b})\}.$$
Let $M$ be an $R$-module. The arithmetic rank
of an ideal ${\frak b}$ of $R$ with respect to $M$, denoted by ${\rm ara}_M({\frak b})$, is defined
the arithmetic rank of the ideal ${\frak b}+\Ann_R(M)/ \Ann_R(M)$ in the ring $R/\Ann_R(M)$.

 \begin{prop}
 \label{2.2}
Let $\frak a$ be an ideal of a Noetherian ring $R$ and $M$ an
$R$-module such that $\dim\Supp(M)\leq 1$ and $\Supp(M)\subseteq V(\frak a)$. Then the
following statements are equivalent:

{\rm(i)} $M$ is $\frak a$-weakly cofinite.

{\rm(ii)} The $R$-modules $\Hom_R(R/\frak a,M)$ and $\Ext^1_R(R/\frak a,M)$
are weakly Laskerian.
 \end{prop}

 \proof The conclusion ${\rm(i)} \Longrightarrow {\rm(ii)}$ is obviously true.
In order to prove that ${\rm(ii)} \Longrightarrow {\rm(i)}$, as
 $$\Ass_R\,\Hom_R(R/\frak a,M)=\Ass_R M$$
and  $\Hom_R(R/\frak a,M)$ is weakly Laskerian,
it follows that $\Ass_R M$ is finite. Now, if $\dim\Supp(M)=0$, then $\Ass_R M= \Supp(M)$, and
so $\Supp(M)$ is also finite. Therefore, in view of definition, $M$ is weakly Laskerian,
and so by \cite[Lemma 2.2]{DM2}, $M$ is $\frak a$-weakly cofinite. Consequently, we may assume $\dim\Supp(M)=1$; and we use induction on $$t:={\rm ara}_M(\frak a)={\rm ara}(\frak a+ \Ann_R
(M)/\Ann_R (M))$$ that $M$ is $\frak a$-weakly cofinite. If $t=0$, then it follows
from definition that $\frak a^n\subseteq {\rm Ann}_R(M)$ for some positive
integer $n$,  and so $M=(0:_M{\frak a}^n)$. Therefore the assertion follows
from \cite[Lemma 2.8]{DM2}.  So assume that $t>0$ and the result has been proved
for all $i\leq t-1$. In view of Lemma 2.6 there exist finitely generated submodules $A$
of $\Hom_R(R/{\frak a},M)$  and $B$ of  $\Ext^1_R(R/{\frak a},M)$
such that the set
\begin{center} 
 $\Omega:= \Supp(\Hom_R(R/{\frak a},M)/A) \bigcup \Supp(\Ext^1_R(R/{\frak a},M)/B).$ 
\end{center} 
is finite. Now, let
$$\mathscr{T}=\{\frak p\in \Supp(M) \mid \dim R/\frak p=1\}\setminus\Omega.$$

It is easy to see that $\mathscr{T}\subseteq {\rm Assh}_R M$, and so
$\mathscr{T}$ is finite. (Note that $\Ass_R M$ is finite.)

In addition, as
$\Omega\subseteq \Supp(M),$
it follows that $$ \max\{\dim\Supp(\Hom_R(R/{\frak a},M)/A), \dim\Supp(\Ext^1_R(R/{\frak a},M)/B)\}\leq 1.$$
Therefore, in view of the prime avoidance theorem it is easy to see that, for each ${\frak p}\in \mathscr{T}$ we have 
${\frak p}\not\subseteq \bigcup _{\frak q \in \Omega}{\frak q}.$
Consequently, it is easily yields that
$$(\Hom_R(R/{\frak a},M)/A)_{\frak p}=0=(\Ext^1_R(R/{\frak a},M)/B)_{\frak p}.$$
Whence for each ${\frak p}\in \mathscr{T}$ the
$R_{\frak p}$-module ${\rm Hom}_{R_{\frak p}}(R_{\frak p}/{\frak a}R_{\frak
p},M_{\frak p})$ is finitely generated, by \cite[Ex. 7.7]{Mat}, and $M_{\frak p}$ is an
$\frak aR_{\frak p}$-torsion $R_{\frak p}$-module, with $\Supp(M)_{\frak p}\subseteq V({\frak p}R_{\frak p})$, and so it follows
that the $R_{\frak p}$-module ${\rm Hom}_{R_{\frak p}}(R_{\frak
p}/{\frak a}R_{\frak p},M_{\frak p})$ is Artinian. Consequently, according
to  Melkersson's results \cite[Theorem 1.3]{Me1} and
\cite[Proposition 4.3]{Me2}, $M_{\frak p}$ is an Artinian and
$\frak aR_{\frak p}$-cofinite $R_{\frak p}$-module. Next, let $\mathscr{T}= \{\frak
p_{1},\dots,\frak p_{n}\}.$ Then by Lemma 3.1(i), we have
$$V({\frak a}R_{\frak p_{j}}) \cap \Att _{R_{\frak p_j}} (M_{{\frak p}_j}) \subseteq V(\frak p_{j}R_{\frak
p_{j}}),$$ for all $j= 1,2, \dots, n$. Next, set

\begin{center}
$\mathscr{U}:= \bigcup_{j=1}^{n} \{\frak q\in \Spec R\mid \frak {q} R_{\frak
p_{j}} \in \Att_{R_{ \frak{p}_{j}}}(M_{{\frak p}_j})\}.$
\end{center}

It is easy to check  that $\mathscr{U} \cap V({\frak a}) \subseteq \mathscr{T}.$\\

On the other hand, since $t={\rm ara}_M({\frak a})\geq 1$, there
exist elements $y_1,  \dots, y_t\in {\frak a}$ such that $${\rm Rad}({\frak a}+{\rm
Ann}_R(M)/{\rm Ann}_R(M))={\rm Rad}((y_1, \dots, y_t)+{\rm
Ann}_R(M)/{\rm Ann}_R(M)).$$  Now, as
%\begin{center}
${\frak a} \not \subseteq \bigcup_{\frak q\in \mathscr{U}\setminus V({\frak a})} \frak q,$
%\end{center}
it follows that
%\begin{center}
$(y_1,\dots, y_t)+{\rm Ann}_R(M) \not \subseteq \bigcup_{\frak q\in
\mathscr{U}\setminus V({\frak a})} \frak q.$
%\end{center}

Furthermore,  for each ${\frak q}\in \mathscr{U}$ we have $\frak {q} R_{\frak
p_{j}} \in \Att_{R_{ \frak{p}_{j}}}(M_{{\frak p}_j}),$ for some
integer $1 \leq j \leq n$. Whence
$${\rm Ann}_R(M)R_{\frak p_{j}}\subseteq {\rm Ann}_{R_{
\frak{p}_{j}}}(M_{{\frak p}_j})\subseteq \frak {q} R_{\frak
p_{j}}.$$  Since $\frak q$ is prime we get that $\Ann_R( M)\subseteq {\frak q}$.
Consequently, it follows from
\begin{center}
$\Ann_R( M) \subseteq \bigcap_{\frak q\in \mathscr{U}\setminus V({\frak a})} \frak
q,$
\end{center}
 that
%\begin{center}
$(y_1, \dots, y_t)\not \subseteq \bigcup_{\frak q\in \mathscr{U}\setminus V({\frak a})}
\frak q.$
%\end{center}
Therefore, by \cite [Ex. 16.8]{Mat} there is $a\in (y_2,\dots,y_t)$
such that
%\begin{center}
$y_1+a\not\in\bigcup_{\frak q\in \mathscr{U}\setminus V({\frak a})} \frak q.$
%\end{center}
Let $x:=y_1+a$. Then $x\in {\frak a}$ and $${\rm Rad}({\frak a}+{\rm Ann}_R(M)/{\rm
Ann}_R(M))={\rm Rad}((x,y_2,...,y_t)+{\rm Ann}_R(M)/{\rm
Ann}_R(M)).$$ Next, let $N:=(0:_Mx)$. Then, it is easy to see that
$${\rm ara}_N({\frak a})={\rm ara}({\frak a}+{\rm Ann}_R(N)/{\rm Ann}_R(N))\leq t-1.$$
(note that $x\in \Ann_RN$), and hence $${\rm Rad}({\frak a}+{\rm
Ann}_R(N)/{\rm Ann}_R(N))={\rm Rad}((y_2, \dots, y_t)+{\rm
Ann}_R(N)/{\rm Ann}_R(N))).$$
 Now, the exact sequence
$$0\longrightarrow N  \longrightarrow M \longrightarrow xM \longrightarrow 0, \,\,\,\,\,\,\,\,\,\,\,\,\,\,\,\,\,\,\,\,\,\,\,\,\,\,\,\,\,\,\,\,\,\,\,\,\,\,\,\,\,\,\,\,\,\,\,\,\,\,\,  (\dag)$$ 
induces an  exact sequence

 $$0\longrightarrow \Hom_R(R/{\frak a}, N) \longrightarrow \Hom_R(R/{\frak a},M)
\longrightarrow \Hom_R(R/\frak a,xM)$$$$\longrightarrow  {\rm
Ext}^{1}_R(R/{\frak a},N) \longrightarrow {\rm
Ext}^{1}_R(R/{\frak a},M),\,\,\,\,\,\,\,\,\,\,\,\,\,\,\,\,\,\,\,\,\,\,\,\,\,\,\,\,\,\,\,\,\,\,\,\,\,\,\,\,\,\,\,\,\,\,\,\,\,\,\,  $$\\
which implies that the $R$-modules ${\rm
Hom}_R(R/{\frak a},N)$ and ${\rm Ext}^1_R(R/{\frak a},N)$ are weakly Laskerian.
Consequently, by the inductive hypothesis, the $R$-module $N$ is
${\frak a}$-weakly cofinite.

 Moreover, the exact sequence $(\dag)$
 induces the exact sequence $$ {\rm
Ext}^{1}_R(R/{\frak a},M) \longrightarrow {\rm
Ext}^{1}_R(R/{\frak a},xM)\longrightarrow {\rm Ext}^{2}_R(R/{\frak a},N),$$ which
implies that the $R$-module ${\rm Ext}^{1}_R(R/{\frak a},xM)$ is weakly Laskerian.\\

 Also, from the exact sequence
$$0\longrightarrow xM  \longrightarrow M \longrightarrow
M/xM \longrightarrow 0$$ we get the exact sequence
$${\rm Hom}_R(R/\frak a,M)
\longrightarrow {\rm Hom}_R(R/\frak a,M/xM)\longrightarrow  {\rm
Ext}^{1}_R(R/\frak a,xM)$$ which implies that the $R$-module ${\rm
Hom}_R(R/\frak a,M/xM)$ is weakly Laskerian.\\
Now, from Lemma 3.1(ii), it is easy to see that the $R_{{\frak p}_j}$-module
$(M/xM)_{\frak p_{j}}$ has finite length for all $j=1,\dots,n$.
Therefore there exists a finitely generated submodule $L_{j}$ of
$M/xM$ such that $$(M/xM)_{\frak p_{j}}= (L_{j})_{\frak p_{j}}.$$ Let
$L:= L_{1}+\cdots+ L_{n}$. Then $L$ is a finitely generated
submodule of $M/xM$ such that

\begin{center}
$\Supp_{R} (M/xM)/L \subseteq {\rm Supp}(M)\setminus \{\frak p_{1},\dots,\frak p_{n}\} \subseteq (\Ass_RM\bigcap\Max(R)) \bigcup \Omega.$
\end{center}

The exact sequence
$$0 \longrightarrow L \longrightarrow M/xM
\longrightarrow (M/xM)/L \longrightarrow 0,$$ provides the following
exact sequence,
$${\rm Hom}_R(R/\frak a,M/xM)
\longrightarrow {\rm Hom}_R(R/\frak a,(M/xM)/L)\longrightarrow  {\rm
Ext}^{1}_R(R/\frak a,L);$$ which implies that 
 $\Hom_R(R/\frak a,(M/xM)/L)$ is weakly Laskerian. 

We now show that $M/xM$ is
a weakly Laskerian $R$-module. To do this, since
the sets  $\Ass_R M\cap \Max(R)$ and $\Omega$  
 are finite,
it follows that the set $\Supp(M/xM)/L$ is finite too. Thus, as $L$ is finitely generated, it follows
from Lemma 2.6 that $M/xM$ is a weakly Laskerian $R$-module. Thus in view of \cite[Lemma 2.6]{DM2} the $R$-module
$M/xM$ is a $\frak a$-weakly cofinite. Now, since the $R$-modules $N=(0:_ M x)$ and $M/xM$ are
$\frak a$-weakly cofinite, it follows from \cite[Lemma 3.1]{Me2} and \cite[Lemma 2.2]{DM2}
that $M$ is $\frak a$-weakly cofinite module. This completes the inductive step. \qed\\

The first application of Proposition 3.2  gives us a characterization of the $\frak a$-weakly cofiniteness of an $R$-module $M$ in terms of the $\frak a$-weakly cofiniteness of the local cohomology modules $H^i_{\frak a}(M)$.
\begin{cor}
 \label{2.4}
Let $R$ be a Noetherian ring, $M$ an  $R$-module and $\frak a$ a one-dimensional ideal of $R$.  Then the following conditions are equivalent:

{\rm(i)} $\Ext^i_R(R/\frak a,M)$ is weakly Laskerian for all $i\leq {\rm cd}(\frak a,M)+1.$

{\rm(ii)} $H^i_{\frak a}(M)$ is  $\frak a$-weakly cofinite for all $i$.

{\rm(iii)} $\Ext^i_R(R/\frak a,M)$ is  weakly Laskerian for all $i$.

{\rm(iv)} $\Ext^i_R(N,M)$ is  weakly Laskerian for all $i\leq {\rm cd}(\frak a,M)+1$ and for any finitely generated
$R$-module $N$ with $\Supp(N) \subseteq V(\frak a)$.

{\rm(v)} $\Ext^i_R(N,M)$ is  weakly Laskerian for all $i\leq {\rm cd}(\frak a,M)+1$ and for some finitely generated
$R$-module $N$ with $\Supp(N) = V(\frak a)$.

{\rm(vi)} $\Ext^i_R(N,M)$ is  weakly Laskerian for all $i$ and for any finitely generated
$R$-module $N$ with $\Supp(N) \subseteq V(\frak a)$.

{\rm(vii)} $\Ext^i_R(N,M)$ is  weakly Laskerian for all $i$ and for some finitely generated
$R$-module $N$ with $\Supp(N) = V(\frak a)$.
\end{cor}
\proof  By a slight modification of the proof of Theroem 2.4, the result follows easily from Proposition 3.2 and Lemmas 2.1, 2.2, by applying  \cite[Lemmas 2.2 and 2.8]{DM2}.   \qed 

\begin{cor}
Let $R$ be a Noetherian ring and let $\frak a, \frak b$ be ideals of $R$ such that
$\frak b\subseteq \Rad(\frak a)$. Let $M$ be a $\frak b$-weakly cofinite $R$-module.

${\rm (i)}$ If $\dim R/\frak a=1$, then the $R$-module $H^i_{\frak a}(M)$ is $\frak a$-weakly cofinite for all $i$.

${\rm (ii)}$ If $\dim R/\frak b=1$, then the $R$-module $H^i_{\frak b}(M)$ is $\frak a$-weakly cofinite for all $i$.
\end{cor}

\proof In order to show that (i), since $\frak b\subseteq \Rad(\frak a)$, it follows that $\Supp(R/\frak a)\subseteq \Supp(R/\frak b)$. On the other hand,
since $M$ is $\frak b$-weakly cofinite it follows from \cite[Lemma 2.8]{DM2}
that $M$ is also $\frak a$-weakly cofinite.  Now since $\dim R/\frak a=1$,  the result follows from Corollary 3.3.

To prove (ii), since $\dim R/\frak b=1$ and $M$ is $\frak b$-weakly cofinite it follows from Corollary 3.3 that
$H^i_{\frak b}(M)$ is $\frak b$-weakly cofinite for all $i$. Now as $\Supp(R/\frak a)\subseteq \Supp(R/\frak b)$ it
follows from \cite[Lemma 2.8]{DM2} that $H^i_{\frak b}(M)$ is $\frak a$-weakly cofinite for all $i$. \qed\\

We are now in a position to use Proposition 3.2 to produce a proof of the main theorem of this section,
which shows that  $\mathscr{C}^1(R, \frak a)_{wcof}$, the category of $\frak a$-weakly cofinite $R$-modules $M$ with $\dim\Supp(M)\leq 1$,  is a full Abelian subcategory of the category of modules.

\begin{thm}
 \label{2.2}
Let $\frak a$ be an ideal of a Noetherian ring $R$. Let
$\mathscr{C}^1(R, \frak a)_{wcof}$ denote the category of $\frak a$-weakly cofinite
$R$-modules $M$ with $\dim\Supp(M)\leq 1$. Then
$\mathscr{C}^1(R, \frak a)_{wcof}$ is an Abelian category.
\end{thm}
\proof Let $M,N\in \mathscr{C}^1(R, \frak a)_{wcof}$ and let $f:M\longrightarrow N$ be
an $R$-homomorphism.  We show that the $R$-modules
$\ker f$ and ${\rm coker} f$ are $\frak a$-weakly cofinite. To this end, the exact sequence
$$0\longrightarrow \ker f \longrightarrow M \longrightarrow
{\rm im} f \longrightarrow 0,$$ induces an  exact sequence
$$0\longrightarrow {\rm Hom}_R(R/\frak a,\ker f) \longrightarrow {\rm
Hom}_R(R/\frak a,M) \longrightarrow {\rm
Hom}_R(R/\frak a,{\rm im} f)$$$$\longrightarrow {\rm Ext}^{1}_R(R/\frak a,\ker f)
\longrightarrow \Ext^{1}_R(R/\frak a,M),$$ that implies the
$R$-modules
\begin{center}
 $\Hom_R(R/\frak a,\ker f)$ and $\Ext^1_R(R/\frak a,\ker f),$
 \end{center}
  are weakly cofinite. Therefore it follows
from Proposition 3.2 that $\ker f$ is $\frak a$-weakly cofinite. Now, by using the exact sequences
$$0\longrightarrow \ker f \longrightarrow M \longrightarrow
{\rm im}f \longrightarrow 0,$$and
$$0\longrightarrow {\rm im}f \longrightarrow N \longrightarrow
{\rm coker}f \longrightarrow 0,$$
we see that ${\rm coker} f$ is also $\frak a$-weakly cofinite, as required. \qed \\

As an immediate consequence of Theorem 3.5, we derive the weakly cofiniteness version of Delfino-Marley's result in \cite{DM} and Kawasaki's result in \cite{Ka2}, which shows that the category of modules weakly cofinite, with respect to an ideal of dimension one in a Noetherian ring, is a full Abelian subcategory of the category of modules.
Following, we let $\mathscr{C}(R, \frak a)_{wcof}$ denote the category of modules weakly cofinite with respect to $\frak a$.
\begin{cor}
Let $\frak a$ be an ideal of a  Noetherian ring $R$ of dimension one.
Then $\mathscr{C}(R, \frak a)_{wcof}$  forms
an Abelian subcategory of the category of all $R$-modules.
\end{cor}
\proof As $\Supp(M) \subseteq \Supp(R/\frak a)$ for all $M\in \mathscr{C}(R, \frak a)_{wcof}$, and $\dim R/\frak a=1$, it follows that $\dim\Supp(M)\leq 1.$  Now the assertion follows from Theorem 3.5. \qed\\
\begin{cor}
 \label{2.2}
Let $\frak a$ be an ideal of a  Noetherian ring $R$ of dimension one.
Let $$X^\bullet:
\cdots\longrightarrow X^i \stackrel{f^i} \longrightarrow X^{i+1}
\stackrel{f^{i+1}} \longrightarrow X^{i+2}\longrightarrow \cdots,$$ be a
complex such that  $X^i\in\mathscr{C}(R, \frak a)_{wcof}$ for all $i\in\Bbb{Z}$.
Then the $i$-th homology module $H^i(X^\bullet)$ is in
$\mathscr{C}(R, \frak a)_{wcof}$.
\end{cor}
\proof The assertion follows from  Corollary 3.6. \qed
\begin{cor}
 \label{2.2}
Let $\frak a=(x_1, \dots, x_n)$ be an ideal of a  Noetherian ring $R$.  Let  $M$ and $N$ be two $R$-modules
such that $N$ is finitely generated and $M$ is $\frak a$-weakly cofinite with $\dim\Supp(M)\leq1$.
Then the $R$-modules $\Ext^i_R(N, M)$, ${\rm Tor}_i^R(N, M)$ and; the Koszul homology module  $H_i(x_1, \dots, x_n; M)$
are $\frak a$-weakly cofinite for all $i$.
\end{cor}
\proof  By considering a finite free resolution $\mathcal{F}_{\bullet}\longrightarrow N$ of $N$, and applying Theorem 3.5 to the complexes 
\begin{center}
$\Hom(\mathcal{F}_{\bullet}, M)$,\,\,\,\,\,\ $\mathcal{F}_{\bullet}\otimes_R M$,\,\,\,\,\,\, ${\bf K}_{\bullet}(x_1, \dots, x_n; M),$
\end{center}  
the assertion follows. \qed\\

We end the paper with the following question:\\

{\bf Question.} Let $\frak a$ be an ideal of a Noetherian ring $R$ and $M$ an
$R$-module such that $\dim\Supp(M)\leq 1$ and $\Supp(M)\subseteq V(\frak a)$. Let $\mathcal{S}$ be a Serre subcategory of the category of $R$-modules. Is the following statements are equivalent ?

{\rm(i)} The $R$-modules  $\Ext^i_R(R/\frak a,M)$ are in $\mathcal{S}$, for all $i\geq 0$.

{\rm(ii)} The $R$-modules $\Hom_R(R/\frak a,M)$ and $\Ext^1_R(R/\frak a,M)$
are in $\mathcal{S}$. 

\begin{center}
{\bf Acknowledgments}
\end{center}

The authors would like to thank from School of Mathematics, Institute for Research in Fundamental Sciences (IPM) for the  financial support. The second author is grateful to the hospitality and facilities offered of the Max-Planck Institut f\"ur Mathematik (Bonn) during the preparation of this paper. 

% ----------------------------------------------------------------

\end{document}